\pdfoutput=1
\RequirePackage{ifpdf}
\ifpdf 
\documentclass[pdftex]{sigma}
\else
\documentclass{sigma}
\fi

\numberwithin{theorem}{section}
\numberwithin{proposition}{section}
\numberwithin{lemma}{section}
\numberwithin{corollary}{section}
\numberwithin{definition}{section}
\numberwithin{example}{section}
\numberwithin{remark}{section}
\numberwithin{note}{section}

\newtheorem{Theorem}{Theorem}[section]
\newtheorem{Corollary}[Theorem]{Corollary}

\newtheorem{Proposition}[Theorem]{Proposition}
{\theoremstyle{definition}

\newtheorem{Example}[Theorem]{Example}
\newtheorem{Remark}[Theorem]{Remark}
\newtheorem{Exercise}[Theorem]{Exercise}}

\newcommand{\coker}{\mathop{\rm coker}\nolimits}

\newcommand{\cha}{\mathop{\rm char}\nolimits}

\newcommand{\supp}{\mathop{\rm Supp}\nolimits}
\renewcommand{\phi}{\varphi}

\newcommand{\rk}{\mathop{\rm rk}\nolimits}

\newcommand{\SL}{{\rm SL}}
\newcommand{\GL}{{\rm GL}}
\newcommand{\SU}{{\rm SU}}

\newcommand{\Hilb}{{\rm Hilb}}

\newcommand{\HM}{\mathop{\rm HM}\nolimits}
\newcommand{\MHM}{\mathop{\rm MHM}\nolimits}
\newcommand{\tMHM}{\widetilde{\mathop{\rm MHM}\nolimits}}
\newcommand{\Perv}{\mathop{\rm Perv}\nolimits}
\newcommand{\Tr}{{\rm Tr}}
\newcommand{\Gr}{\mathop{\rm Gr}\nolimits}

\newcommand{\Sym}{{\rm Sym}}

\renewcommand{\tilde}{\widetilde}

\newcommand{\F}{\mathcal F}
\newcommand{\hH}{\mathbb H}
\newcommand{\cH}{\mathcal H}
\renewcommand{\P}{\mathbb P}

\newcommand{\M}{\mathcal M}
\newcommand{\oM}{\overline{\mathcal M}}

\newcommand{\D}{\mathcal D}

\newcommand{\A}{{\mathbb A}}

\newcommand{\E}{\mathcal E}
\newcommand{\R}{\mathbb R}
\newcommand{\cR}{\mathcal R}

\newcommand{\cN}{\mathcal N}
\newcommand{\cC}{\mathcal C}
\newcommand{\cL}{\mathcal L}
\newcommand{\C}{\mathbb C}

\newcommand{\Z}{\mathbb Z}

\renewcommand{\O}{\mathcal O}
\renewcommand{\L}{\mathbb L}

\newcommand{\Q}{\mathbb Q}

\renewcommand{\sl}{{\mathfrak sl}}

\newcommand{\half}{\frac{1}{2}}

\begin{document}

\allowdisplaybreaks

\renewcommand{\thefootnote}{$\star$}

\renewcommand{\PaperNumber}{088}

\FirstPageHeading

\ShortArticleName{Nekrasov's Partition Function and Ref\/ined DT Theory}

\ArticleName{Nekrasov's Partition Function and Ref\/ined\\ Donaldson--Thomas Theory: the Rank One Case\footnote{This
paper is a contribution to the Special Issue ``Mirror Symmetry and Related Topics''. The full collection is available at \href{http://www.emis.de/journals/SIGMA/mirror_symmetry.html}{http://www.emis.de/journals/SIGMA/mirror\_{}symmetry.html}}}

\Author{Bal\'azs SZENDR\H{O}I}

\AuthorNameForHeading{B.~Szendr\H{o}i}

\Address{Mathematical Institute, University of Oxford, UK}
\Email{\href{mailto:szendroi@maths.ox.ac.uk}{szendroi@maths.ox.ac.uk}}
\URLaddress{\url{http://people.maths.ox.ac.uk/szendroi/}}

\ArticleDates{Received June 12, 2012, in f\/inal form November 05, 2012; Published online November 17, 2012}

\Abstract{This paper studies geometric engineering, in the simplest possible case of rank one
(Abelian) gauge theory on the af\/f\/ine plane and the resolved conifold. We recall the identif\/ication between
Nekrasov's partition function and a version of ref\/ined Donaldson--Thomas theory, and study the
relationship between the underlying vector spaces. Using a~purity result, we identify the vector
space underlying ref\/ined Donaldson--Thomas theory on the conifold geometry as the exterior space of the space of
polynomial functions on the af\/f\/ine plane, with the (Lefschetz) $\SL(2)$-action on the threefold side
being dual to the geometric $\SL(2)$-action on the af\/f\/ine plane.
We suggest that the exterior space should be a module for the
(explicitly not yet known) cohomological Hall algebra (algebra of BPS states) of the conifold.}


\Classification{14J32}

\renewcommand{\thefootnote}{\arabic{footnote}}
\setcounter{footnote}{0}

\section{Introduction}

We study an instance of geometric engineering~\cite{KKV}. Our starting point is Type II string theory on
a real 10-dimensional spacetime $X\times \C^2$, the product of the f\/lat Calabi--Yau surface $\C^2$ and
a~local Calabi--Yau threefold $X$. This threefold is determined by a f\/inite subgroup $\Gamma<\SU(2)$;
it is the minimal resolution $X\to \O_{\P^1}(-1, -1)/\Gamma$
of the singular threefold $\O_{\P^1}(-1, -1)/\Gamma$, where $\Gamma$ acts f\/iberwise
on the resolved conifold $\O_{\P^1}(-1, -1)$, f\/ixing the zero-section.
Integrating out the $X$-direction leads to supersymmetric gauge theory on~$\C^2$ with gauge group
$G=G(\Gamma)$, the simple group of type~A,~D or~E corresponding to the type of~$\Gamma$. On
the other hand, integrating out the $\C^2$-direction gives a $\sigma$-model, or equivalently
a version of $U(1)$ gauge theory on $X$.

The aim of this paper is to revisit the construction of the partition functions $Z_{\C^2}$ and
$Z_X$ on the two sides and their identif\/ication, mainly in the simplest possible case when $\Gamma$ is the trivial
group. As it turns out, both partition functions are characters
of representations of tori (Hilbert series of graded vector spaces). The underlying vector spaces are
the symmetric, respectively the exterior space of the space of functions on $\C^2$. On the four-dimensional
side this is not new, but it is a surprising fact on the six-dimensional side, and suggests that
a version of Koszul duality might lurk behind the mathematics of geometric engineering.
It would be very interesting to see more examples in action to substantiate this claim;
we comment on the dif\/f\/iculties below.
Slightly more concretely, we suggest at the end of the paper that the alternating space
$\Lambda^*\C[x,y]$ should be a module for the Kontsevich--Soibelman CoHA (cohomological Hall algebra)
$\cH(Q,W)$ attached to the conifold quiver $(Q,W)$; an explicit description of this algebra is currently
lacking but would be desirable.

\section{The four-dimensional partition function}

While it would be possible to remain more general for a while at least, let us make a simplifying
assumption right away and assume that $\Gamma$ is the cyclic group of order $r$, embedded diagonally into $\SU(2)$.
Thus, on the four-dimensional side we are considering $\SU(r)$ gauge theory on the complex plane $\C^2$.
This is certainly well-def\/ined for $r>1$; for $r=1$ we recall the interpretation below.

The partition function of the theory is a sum of integrals over a collection of moduli spaces, the spaces
of $\SU(r)$ instantons on $\R^4$ of various charges $k$. It is well known that, after framing the instantons,
their moduli space is dif\/feomorphic to the moduli spaces
\[\M^\circ_{r,k} =\left\{ \begin{array}{c|c}(\E, \phi) & \begin{array}{l} \E \mbox{ vector bundle on }\P^2 \\ \rk\E=r, \ c_2(\E)=k, \ \phi\colon \E|_{l_\infty}\cong \O_{l_\infty}^{\oplus r}\end{array}\end{array}\right\}\Big/\sim \]
of framed rank-$r$ bundles $\E$ on $\P^2$ of charge $k$.
Here $l_\infty\cong\P^1$ is the complement of $\C^2$ in $\P^2$.
Note that the parametrized bundles $\E$ are indeed $\SL(r)$-bundles, since the framing
isomorphism automatically trivializes the determinant. This space admits two partial
compactif\/ications $\M_{r,k}$ and $\oM_{r,k}$. First, we have the space
\[\M_{r,k} =\left\{  \begin{array}{c|c}(\E, \phi) & \begin{array}{l} \E \mbox{ torsion-free sheaf on }\P^2 \\  \rk\E=r, \ c_2(\E)=k, \ \phi\colon \E|_{l_\infty}\cong \O_{l_\infty}^{\oplus r}\end{array}\end{array}\right\}\Big/\sim, \]
the moduli space of framed rank-$r$ torsion-free sheaves
on $\P^2$ of charge $k$; this is the analogue of the Gieseker compactif\/ication of the moduli of bundles
for a projective surface. $\M_{r,k}$~is a~nonsingular holomorphic sympletic variety of dimension~$2kr$.
The second space $\oM_{r,k}$ is the analogue of the Uhlenbeck compactif\/ication,
and can be constructed as an af\/f\/ine GIT quotient; for details, see~\cite[Section 2]{NY1}.

The four-dimensional gauge theoretic partition function (for pure gauge theory, in the absence of matter) is
\begin{equation}\label{eg_Z_def}
Z_{\C^2,r}(\Lambda) = \sum_{k\geq 0} \Lambda^k \int_{\M_{r,k}} 1,
\end{equation}
the generating function of symplectic volumes of the (non-singular) moduli spaces $\M_{r,k}$.
This is an ill-def\/ined expression, since $1$ is not a top-dimensional form, and integration happens over
noncompact spaces~$\M_{r,k}$.

As Nekrasov~\cite{N} discovered, one can ``renormalize''
both these problems by considering equivariant integration with respect to
the torus $T=(\C^*)^2\times (\C^*)^{r-1}$. Here the f\/irst factor acts via its geometric action
on $\C^2$, which extends to an action on $\M_{r,k}$. The second component $(\C^*)^{r-1}$, the maximal torus
of $\SL(r)$, acts on the framing $\phi$, f\/ixing the trivialization of the determinant.

We will in fact be interested in a K-theoretic version\footnote{What I call $Z_{\C^2,r}$ is in fact called the 5-dimensional
partition function in the physics literature, since it arises from studying M-theory on a circle bundle over $X\times \C^2$.
I will abuse language and will continue to call it the four-dimensional partition function, since it is still naturally
associated to the complex surface.}
of the partition function, also introduced
in \cite{N} and studied in detail in \cite{NY1, NY2}.
In the K-theoretic version, the integrand $1$ gets
interpreted as the unit K-theory class, the class of the structure sheaf $\O$; integration gets replaced
by pushforward to the point. Thus, the partition function computes the generating series of equivariant coherent
cohomologies of $\O$. This sort of generating series had been considered earlier
in a related context under the name of four-dimensional Verlinde formula in~\cite{LMNS}.

Introducing a basis $q_1,q_2, a_1, \ldots, a_{r-1}$ for the space of $T$-characters,
the $K$-theoretic partition function is
\begin{equation} Z_{\C^2,r}(q_i, a_j,\Lambda) = \sum_{k\geq 0} \Lambda^k \cha_T R(\pi_{r,k})_* \O_{\M_{r,k}} \in \Z(q_i, a_j)[[\Lambda]].
\label{eq_nek_def}
\end{equation}
Here $\pi_{r,k}\colon  \M_{r,k}\to \{*\}$ is the structure morphism, $\cha_T V$ denotes the $T$-character
of a representation $V$ of the torus~$T$, and $\cha_T R(\pi_{r,k})_*$ is shorthand for
$\sum_i (-1)^i \cha_T R^i(\pi_{r,k})_*$, the torus character of an (a priori) virtual representation of $T$.
The formula~(\ref{eq_nek_def}) gives a well-def\/ined expression, since it is easy to see
that each $T$-weight space of each $R^i(\pi_{r,k})_* \O_{\M_{r,k}}$ is f\/inite
dimensional~\cite[Section 4]{NY1}. It can be shown that the answer is indeed a rational function of
the variables $q_i$, $a_j$.

It is known that $\pi_{r,k}$ factors through a map $\sigma_{r,k}\colon \M_{r,k}\to \oM_{r,k}$ to the Uhlenbeck-type
space. Pushforward under $\sigma_{r,k}$ gives no higher cohomology~\cite[Lemma 3.1]{NY1};
also $\oM_{r,k}$ is af\/f\/ine, so there is no higher cohomology there either. Hence the partition
function is simply
\begin{equation} Z_{\C^2,r}(q_i, a_j,\Lambda) = \sum_{k\geq 0} \Lambda^k \cha_T H^0\left( \oM_{r,k}, \O_{\oM_{r,k}}\right) \in \Z(q_i, a_j)[[\Lambda]]. \label{eq_nek_H0}
\end{equation}

Let us restrict further to the case $r=1$, which in any case requires some further explanation.
In this case we have ``$\SU(1)$ gauge theory'' on $Y=\C^2$. This gets interpreted as
the theory of framed rank-one sheaves with trivial determinant on $\P^2$. Since there is only one line bundle
on~$\P^2$ with trivial determinant, the moduli space $\M_{1,k}$ of framed rank-one
torsion-free sheaves on~$\P^2$ can be identif\/ied with $\Hilb^k(\C^2)$, the
Hilbert scheme of~$k$ points on~$\C^2$. The corresponding Uhlenbeck space $\oM_{1,k}$ is the symmetric product $S^k(\C^2)$;
$\sigma_{1,k}$ is the Hilbert--Chow morphism.
There are no~$a_j$ parameters. The partition function in this case can be computed explicitly in closed form.

\begin{Proposition}[\cite{NY1, NO}] \label{prop_U1} We have
\begin{equation}
Z_{\C^2, r=1}(q_1, q_2,\Lambda)= \prod_{i_1, i_2\geq 0}\big(1-q_1^{i_1}q_2^{i_2}\Lambda\big)^{-1}.\label{eq_nekU1}
\end{equation}
\end{Proposition}
\begin{proof}
\begin{gather*}
Z_{\C^2, r=1}(q_1, q_2,\Lambda)   =   \sum_{k\geq 0} \Lambda^k \cha_T H^0(\O_{S^k(\C^2)})
  =  \sum_{k\geq 0} \Lambda^k \cha_T \C[x_1, \ldots, x_k, y_1, \ldots, y_k]^{S_k} \\
\phantom{Z_{\C^2, r=1}(q_1, q_2,\Lambda) }{}
  =   \sum_{k\geq 0} \Lambda^k \cha_T  S^k \C[x, y]  =  \cha_{T\times\C^*} S^* \C[x, y]
 =  \prod_{i_1, i_2\geq 0}(1-q_1^{i_1}q_2^{i_2}\Lambda)^{-1}.
\end{gather*}
Here in the penultimate line, $S^* \C[x, y]$ is treated as a triply-graded space, graded by $x$-weight, $y$-weight
and polynomial weight with respect to the outer symmetric power operation. A triple grading corresponds to
an action of a rank-three torus $T\times\C^*$, and the character is the Hilbert series.
\end{proof}

In the higher rank case, there is no known closed formula for $Z_{\C^2,r}(q_i, a_j,\Lambda)$. Through torus
localization, it can be computed as a sum over the $T$-f\/ixed points of the spaces $\M_{r,k}$, parametrized by
$r$-tuples of partitions~\cite{N, NO}. It is shown in~\cite{NY2} that it satisf\/ies a system of functional equations called
the blowup equations, whose solution is unique.

\section{The six-dimensional partition function}

Under geometric engineering, the four-dimensional partition function~$Z_{\C^2,r}$ should correspond to a partition function
built out of invariants of the Calabi--Yau threefold $X$. The gauge group on the threefold $X$
is always going to be $\SU(1)$, in other words we will be looking at a version of rank one sheaf theory
with trivial determinant. To f\/ind out precisely which version, let us restrict to the case $r=1$ again.

From the four-dimensional theory we obtain the partition function $Z_{\C^2, r=1}(\Lambda, q_1, q_2)$.
Taking its inverse and specializing, consider
\[ Z(q,T) = Z_{\C^2, r=1}\left(\Lambda= -Tq, q_1=q_2=-q\right)^{-1} = \prod_{n\geq 1} (1-(-q)^nT)^n.
\]
This is a very well-known expression, the reduced topological string partition function~\cite{GV}
of the resolved conifold $X=\O_{\P^1}(-1, -1)$. After a further change of variables~\cite{MNOP},
the function $\log Z(q=-e^{i\hbar}, T)$ gives the full Gromow--Witten potential of the resolved conifold $X$,
with $\hbar$ being the genus parameter.

On the other hand, the precise geometric interpretation of the coef\/f\/icients $Z(q,T)$ is given by
\begin{Theorem}[\cite{NN}] We have
\[ \prod_{n\geq 1} (1-(-q)^n T)^n = \sum_{l,m\geq 0} P_{l,m} T^l q^m,
\]
where $P_{l,m}$ are the pairs invariants or Pandharipande--Thomas $($PT$)$ invariants~{\rm \cite{PT}}
of the resolved conifold~$X$.
\end{Theorem}

Let us recall from~\cite{PT} how the invariants $P_{l,m}$ are def\/ined. They are the enumerative invariants
associated to a collection of highly singular moduli spaces $\cN_{l,m}$, which carry a perfect obstruction theory.
These spaces are the moduli spaces of stable pairs
\[\cN_{l,m} = \left\{ \begin{array}{c|c}(\F, s) & \begin{array}{lc} \F \mbox{ pure 1-dimensional sheaf with proper support on }X \\ s\colon \O_X\to \F \mbox{ a section} \\
\supp(\F)=l[\P^1], \ \chi(\F)=m, \  \dim\supp\coker(s)=0 \end{array}\end{array}\right\}\Big/\sim,\]
where the restriction of the cokernel of $s$ having zero-dimensional support is the stability condition.
Note that the sheaf $\F$, having necessarily proper support, must be supported on a multiple of
the zero-section in $X$. As observed by Bridgeland~\cite{B}, the spaces  $\cN_{l,m}$ represent a moduli problem involving perverse
coherent sheaves on $X$. Perverse coherent sheaves are complexes of coherent sheaves, which belong to the heart of a
$t$-structure on the derived category of sheaves on $X$ dif\/ferent from the standard one.

Our aim is to f\/ind a one-parameter ref\/inement of the numbers $P_{l,m}$, in order to obtain a~ref\/inement of the topological string partition function of the conifold which matches the full
Nekrasov's formula. The following result is crucial for further progress.

\begin{Theorem} \label{thm_super}
The moduli spaces $\cN_{l,m}$ are proper. Moreover, they are global degeneracy loci:
there exist smooth varieties $N_{l,m}$
equipped with regular functions $f_{l,m}\colon N_{l,m}\to\C$ so that, scheme-theoretically,
\[ \cN_{l,m} = \left\{ df_{l,m}=0\right\} \subset N_{l,m}
\]
are the degeneracy loci of these functions.
\end{Theorem}
\begin{proof} The f\/irst statement follows from the fact that the only proper curve in~$X=\O_{\P^1}(-1, -1)$ is the zero section,
so the reduced support of a stable pair cannot move.

To prove the
second statement, we need to recall the quiver interpretation of the moduli
spaces attached to the resolved conifold~\cite{Sz}.
First of all, let $Q$ be the conifold quiver, consisting of vertex set $V=\{0,1\}$ and edge set $\{a_{01}, b_{01}, a_{10}, b_{10}\}$,
with edges labelled $ij$ pointing from vertex $i$ to vertex $j$.
Consider the superpotential~\cite{KW} \[W=a_{01}a_{10}b_{01}b_{10} -a_{01}b_{10}b_{01}a_{10}.\]
Let also $\tilde Q$ be the framed (or extended)
quiver with vertex set $\tilde V=\{0,1,\infty\}$, an extra edge~$i_{\infty 0}$ and the same superpotential $\tilde W=W$.
Then by~\cite{NN}, the moduli spaces $\cN_{l,m}$ can be identif\/ied with stable representations of the quiver $\tilde Q$
with relations arising from formal partial derivatives of $\tilde W$ with respect to the various edges, and a dimension vector
${\mathbf d}=(d_0, d_1, d_\infty=1)$ where $(d_0, d_1)$ depends on~$l$,~$m$. Stability is taken with respect to a specif\/ic
stability condition, determined by a particular (limiting) value of a stability parameter $\theta\in \R^V$ (see more on this below).

Now let $N_{l,m}$ be the moduli space of $\theta$-stable representations of the quiver $\tilde Q$ with the same
dimension vector ${\mathbf d}$, but with no relations. Since $\theta$ is chosen generically, this is a smooth
quasiprojective variety. Let also $f_{l,m}\colon N_{l,m}\to \C$
be def\/ined by evaluating the expression $\Tr(W)$ on representations. It is now well known that the scheme-theoretic equations of
$\cN_{l,m}\subset N_{l,m}$ are indeed given by $d\Tr(W)=0$.
\end{proof}

\begin{Remark}\label{rem:lb} There is a further point worth mentioning in connection with this construction. As
well as the smooth GIT quotient $N_{l,m}$, there is also the af\/f\/ine GIT quotient $\overline{N}_{l,m}$ together
with a contraction morphism $N_{l,m}\to \overline{N}_{l,m}$ (an analogue of the map $\sigma_{r,k}$ in the
four-dimensional situation). This is a projective morphism, equipped by
construction with a relatively ample line bundle (coming from the choice of GIT stability).
On the other hand, $\cN_{l,m}$, being proper, must sit in a f\/ibre of this contraction. Thus $\cN_{l,m}$
comes automatically equipped with a polarization, a chosen ample line bundle $\cL_{l,m}$.
\end{Remark}

Theorem~\ref{thm_super} implies that the singular moduli space $\cN_{l,m}$ acquires a topological coef\/f\/icient
system
\[\phi_{l,m} =\phi_{f_{l,m}} \Q_{N_{l,m}}[\dim N_{l,m}] \in \Perv_\Q(\cN_{l,m}),\]
the {\em perverse sheaf of vanishing cycles} of the regular function $f_{n,l}$ on $N_{l,m}$. This
perverse $\Q$-sheaf is well known to live on the (reduced) degeneracy locus. In fact, we wish to use
Hodge theory, so we consider the canonical lift
\[\Phi_{l,m} =\phi^H_{f_{n,l}} \Q^H_{N_{l,m}}(\dim N_{l,m}/2)[\dim N_{l,m}] \in \tMHM_\Q(\cN_{l,m}),\]
the mixed Hodge module of vanishing cycles of $\phi_{f,n}$.
To make sense of the above expression, we need to extend the category of mixed Hodge modules by a half-Tate object;
see Appendix~\ref{app_weights}.

Now consider the (hyper)cohomology $\hH^*(\cN_{l,m}, \Phi_{l,m})$. As the cohomology of a
mixed Hodge module, it carries a weight f\/iltration $W_n$.
Consider the weight polynomial\footnote{We assume here for simplicity of exposition
that the semisimple part of the monodromy
acts trivially on the cohomology, as will be the case in the example where we apply this formalism. In general, the weight f\/iltration needs to be shifted on the part of the cohomology where the semisimple monodromy acts nontrivially.}
 (for details and examples, see Appendix~\ref{app_weights} again)
\begin{equation}\label{Wdef} W\big(\cN_{l,m}, \Phi_{l,m}; t^\half\big) = \sum_{n\in \Z} t^{\frac{n}{2}} \sum_{i\in\Z} (-1)^i \dim \Gr^W_{n} \hH^i(\cN_{l,m}, \Phi_{l,m})\in\Z\big[t^{\pm\half}\big].
\end{equation}
Here $\Gr^W_n$ denotes the $n$-th graded piece under the weight f\/iltration. Note that the weight polynomial
is usually def\/ined on compactly supported cohomology, but this makes no dif\/ference here, as $\cN_{l,m}$ is
proper, and $\Phi_{l,m}$ is self-dual under Verdier duality.

The following result shows that we are indeed considering here a one-parameter ref\/inement of the numerical
invariants introduced earlier.
\begin{Proposition} We have
\[P_{l,m}  = W\big(\cN_{l,m}, \Phi_{l,m}; t^\half=1\big).\]
\end{Proposition}
\begin{proof} Setting $t^\half=1$ in the weight polynomial, we are ignoring the weight f\/iltration,
and thus computing the Euler characteristic of the perverse sheaf of vanishing cycles of the function~$f_{l,m}$ on~$N_{l,m}$. As discussed in~\cite{Be}, this computes the integral (weighted Euler characteristic)
of a~canonical constructible function, the so-called Behrend function, on the degeneracy locus~$\cN_{l,m}$.
On the other hand, the numerical invariant $P_{l,m}$ arises from a symmetric perfect obstruction theory on~$\cN_{l,m}$.
By the main result of~\cite{Be}, the degree of the associated virtual fundamental class
agrees with the Euler characteristic of $\cN_{l,m}$ weighted by the Behrend function.
\end{proof}

Consider the generating series
\[ Z_X\big(q,T,t^\half\big) = \sum_{l,m\in\Z} W\big(\cN_{l,m}, \Phi_{l,m}; t^\half\big) T^l q^m.
\]
This series is computed in the following result.
\begin{Theorem}[\cite{MMNS}] \label{thm_mmns} We have
\begin{equation}
Z_X\big(q,T,t^\half\big) = \prod_{m\geq 1}\prod_{j=0}^{m-1} \big(1- (-1)^m t^{-\frac m2+\half+j} q^m T \big).
\label{eq_mmns}
\end{equation}
\end{Theorem}

\begin{Remark} Our paper~\cite{MMNS} works in a slightly dif\/ferent framework, considering
a ring-valued rather than cohomological invariant, taking values in a version of the motivic ring of varieties,
thus adopting in the approach of~\cite{KS1}. However, there is a homomorphism from the motivic ring
to the polynomial ring $\Z[t^{\pm\half}]$, given by the weight polynomial, since the weight polynomial
respects the fundamental def\/ining relation $[X]=[X\setminus Z]+[Z]$ of the motivic ring,
where $Z\subset X$ is a closed subvariety. Hence the results of~\cite{MMNS} indeed apply.
\end{Remark}

Thus, comparing~\eqref{eq_nekU1} with~\eqref{eq_mmns}, we obtain a full interpretation
of (the inverse of) the Nekrasov partition function in the simplest case $r=1$:
\begin{Corollary}\label{cor_identify}
The four-dimensional and six-dimensional generating series are related by the equality
\begin{equation} Z_{\C^2, r=1}(q_1, q_2,\Lambda) = Z_X\big(q,T,t^\half\big)^{-1}
\label{fullZinverse}
\end{equation}
under the change of variables $q_1=-t^\half q$, $q_2=-t^{-\half} q$, $\Lambda=-Tq$.
\end{Corollary}

\begin{Remark} The search for, and development of, a ref\/ined version of Donaldson--Thomas theory has always
been motivated and informed by the relationship to Nekrasov's partition function~\cite{HIV, IK, IKV}.
The fact that the motivic, or equivalently the weight polynomial version of DT theory should be the
right ref\/inement was suggested independently in~\cite{BBS, DS} and~\cite{DG}.

Ongoing work of Nekrasov and Okounkov~\cite{NO2}, announced in~\cite{O}, gives a K-theoretic interpretation
of the six-dimensional partition function $Z_X$, which, under suitable assumptions, is compatible\footnote{Maulik D., Okounkov A., work in progress.} 
with the Hodge-theoretic interpretation given above.
\end{Remark}

\begin{figure}[ht]
\centering
\includegraphics{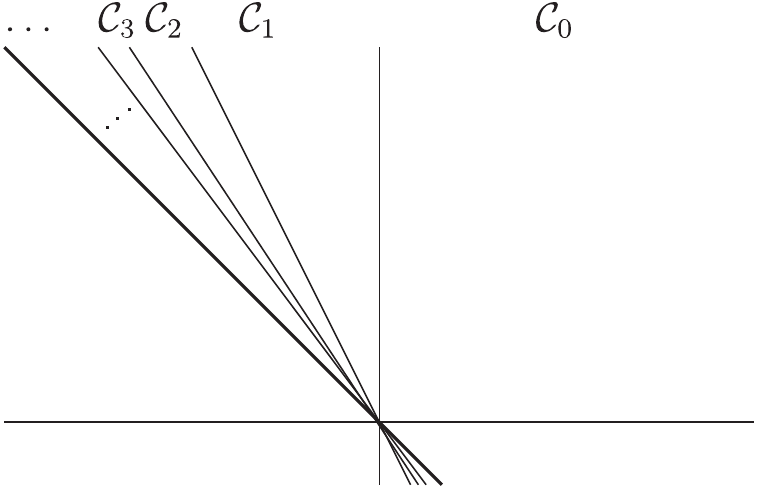}
\caption{Some of the walls and chambers in the space of stability conditions on the framed conifold quiver,
with the thick line representing the D0/D6 wall.}\label{Fig1}
\end{figure}

A variant of Theorem~\ref{thm_mmns} will be of later use. Recall the interpretation of the spaces $\cN_{l,m}$
as spaces of stable representations of the conifold quiver $(\tilde Q, \tilde W)$.
The specif\/ic stability condition depends on a stability parameter $\theta\in \R^2$, and the
spaces $\cN_{l,m}$ arise when $\theta$ takes a value in a~certain limiting position~\cite{NN}. More precisely,
there is a set of open cones $\cC_0, \cC_1, \ldots \subset \R^2$, as in Fig.~\ref{Fig1},
and moduli spaces $\cN^\theta_{l,m}$ for $\theta$
inside any of the open chambers $\cC_n$, such that Theorem~\ref{thm_super} continues to hold.
Also, for each f\/ixed $l$,~$m$, the moduli spaces  $\cN^\theta_{l,m}$ stabilize as $\theta\in\cC_n$ with $n\to\infty$, and
agree with the stable pair moduli spaces $\cN_{l,m}$. Now for any~$\theta$ inside a chamber, we can consider
\[
Z^\theta_X\big(q,T,t^\half\big) = \sum_{l,m\in\Z} W\big(\cN^\theta_{l,m}, \Phi^\theta_{l,m}; t^\half\big) T^l q^m.
\]

\begin{Theorem}[\cite{MMNS}]\label{thm_mmns:theta}
For $n\geq 0$ and $\theta\in\cC_n$, we have
\begin{equation*}
Z^\theta_X\big(q,T,t^\half\big) = \prod_{m=1}^n\prod_{j=0}^{m-1} \big(1- (-1)^m t^{-\frac m2+\half+j} q^m T \big).
\end{equation*}
\end{Theorem}

Beyond the chambers $\cC_n$ as $n\to\infty$, is the D0/D6 wall, the other side of which the moduli spaces
$\cN^\theta_{l,m}$ are no longer compact. Immediately on the other side, (inf\/initesimally) close to the wall,
these moduli spaces are in fact the actual Donaldson--Thomas moduli spaces representing rank-one sheaves with
trivial determinant (and compactly supported quotient)~\cite{NN}. In particular, the Hilbert schemes of points
of the threefold $X$ make an appearance. The partition functions here include further, ref\/ined MacMahon-type
factors~\cite{MMNS}.

Let us return f\/inally to the case $r>1$. First consider the specialization where the torus parameters $q_1, q_2$ are identif\/ied
as above. In this case, it is known~\cite{IK} that the K-theoretic version of the gauge theory partition function
reproduces the full Gromow--Witten or reduced Donaldson--Thomas series of the threefold~$X$. (I am simplifying here:
the threefold $X$ def\/ined in the Introduction is a f\/ibration over $\P^1$ by resolved $A_{r-1}$ surface singularitites.
For $r>1$, there is in fact a family of such threefolds~\cite[Fig.~10]{HIV},
whose Gromow--Witten potentials dif\/fer by the inclusion of certain framing
terms. On the surface side, this dif\/ference is accounted for~\cite{T} by a~change in the integrand in~\eqref{eg_Z_def},
from the trivial class to the class of a power of the determinant line bundle.)

The ref\/ined case is largely open for $r>1$. Even for $r=2$, when $X=\O_{\P^1\times\P^1}(-2,-2)$ is a much-studied
local Calabi--Yau threefold, there appear to be signif\/icant challenges. The ref\/ined topological vertex formalism gives
a combinatorial answer~\cite[Section~5.5]{IKV}, which can be matched with the Nekrasov partition
function~\cite[Section~5.1.1]{IKV}. The very recent paper~\cite{CKK} uses localization techniques to compute some
ref\/ined PT invariants, making some assumptions which are justif\/ied using~\cite{NO2, O}.
On the quiver side however, the wall crossing picture of~\cite{NN}, used extensively in~\cite{MMNS},
is a lot more complicated and hardly understood. When trying to match the quiver picture with the geometry,
there is no clear understanding how f\/ind DT or PT moduli spaces starting from the quiver.
It would be very interesting to f\/ind a way to calculate all ref\/ined PT invariants for this example.

\section{Purity and Hard Lefschetz}

In this section, I will discuss a result which will be used in the f\/inal section. To motivate it, let us assume
f\/irst that $X=\{df=0\}\subset N$ is a global degeneracy locus of a smooth function $f\colon N\to\C$, and let us
moreover assume that $X$ is smooth (and reduced) at every point. Then it follows from the Morse lemma
with parameters that the vanishing
cycle perverse sheaf is just a~local system on $X$ of rank one. Assuming also that $X$ is simply connected,
the coef\/f\/icient system~$\Phi$ associated to $f$ is just a shift of the trivial sheaf (Hodge module) $\Q_X^H$.
Assuming f\/inally that~$X$ is proper, by Deligne's result,
the cohomology $\hH^*(X, \Phi)$ carries a pure Hodge structure, where the weight f\/iltration agrees with
the degree f\/iltration: the weight of a cohomology class equals its degree.
In particular, the weight polynomial $W(X,\Phi)$ is just a shift of the topological Poincar\'e polynomial of $X$.
Moreover, if $L\in H^2(X,\Z)$ is any ample class on $X$, then the maps
\[\begin{array} {@{}rccl} L_k \colon & \hH^{-k}(X, \Phi) & \to & \hH^k(X, \Phi)(k)\\
& \alpha& \mapsto & \alpha\cup L^k
\end{array}
\]
are isomorphisms by the Hard Lefschetz theorem, and can be used to endow $\hH^*(X, \Phi)$ with an $\sl(2)$-action.

Our moduli spaces $\cN_{l,m}$ are not at all smooth, and the coef\/f\/icient systems are nontrivial. However, we have the
following general result.

\begin{Theorem}[\cite{DMSS}]
\label{thm:pure} Let $f\colon X\to\C$ be a regular function on a smooth quasi-projective variety. Assume
that $X$ carries a $\C^*$-action, so that $f$ is equivariant with respect to the weight-$1$ action of $\C^*$ on
the base $\C$. Assume also that the critical locus $Z=\{df=0\}$ is proper, carrying an ample line bundle $\cL$.
\begin{enumerate} \itemsep=0pt
\item[$1.$] The Hodge structure on the cohomology $\hH^i(Z, \Phi_{f})$ is pure of weight~$i$;
equivalently, the weight filtration on $\hH^*(X, \Phi_{f})$ agrees with the degree filtration.
\item[$2.$] The ample class $\cL$ defines, by cup product as above, Hard Lefschetz-type isomorphisms
\[ \cL_k \colon \ \hH^{-k}(Z, \Phi_{f}) \to  \hH^{k}(Z, \Phi_{f})(k)
\]
leading to an $\sl(2)$-action on $\hH^*(\cN_{l,m}, \Phi_{l,m})$.
\end{enumerate}
\end{Theorem}

As a consequence, we obtain

\begin{Corollary}\label{cor:pure} For all $l$, $m$ we have that
\begin{enumerate}\itemsep=0pt
 \item[$1)$] the Hodge structure on the cohomology $\hH^i(\cN_{l,m}, \Phi_{l,m})$ is pure of weight~$i$;
equivalently, the weight filtration on $\hH^*(\cN_{l,m}, \Phi_{l,m})$ agrees with the degree filtration;
\item[$2)$] the ample class $\cL_{l,m}$ defines, by cup product as above, Hard Lefschetz-type isomorphisms
\[ \cL_k \colon \ \hH^{-k}(\cN_{l,m}, \Phi_{l,m}) \to  \hH^{k}(\cN_{l,m}, \Phi_{l,m})(k)
\]
leading to an $\sl(2)$-action on $\hH^*(\cN_{l,m}, \Phi_{l,m})$.
\end{enumerate}
\end{Corollary}

\begin{proof} By Theorem~\ref{thm_super}, the moduli spaces $\cN_{l,m}$ are proper degeneracy loci,
and they carry ample line bundles by Remark~\ref{rem:lb}. To conclude, it is suf\/f\/icient to observe that the Klebanov--Witten
superpotential is linear in each of the variables, so the functions $f_{l,m}$ can be made homogeneous of degree one
by acting by $\C^*$ on the matrix corresponding to either of the four arrows of the conifold quiver.
\end{proof}

\begin{Remark} For a pure Hodge module $\Phi$ on a proper scheme $X$ of some f\/ixed weight, the cohomology is also pure
by the decomposition theorem, and one also has the Lefschetz package~\cite[Theorem~5.3.1, Remark~5.3.12]{S}.
However, the mixed Hodge module $\Phi_{l,m}$ on $\cN_{l,m}$ is known not to be of f\/ixed weight in some cases;
see Example~\ref{ex_N24} below. There is a similar, but dif\/ferent example of a mixed Hodge module with pure cohomology
on a singular projective variety in~\cite[Section~6]{E}.
\end{Remark}

\begin{Example}  Consider f\/irst the case $l=1$, corresponding to pairs invariants with the curve class
having multiplicity $1$. As shown in~\cite[Section~4.1]{PT}, the geometry of the corresponding
moduli spaces $\cN_{1,m}$ is simple: for $m\geq 1$,
$\cN_{1,m}$ it is the space of nonzero sections, up to scale, of $\O_{\P^1}(m-1)$. Thus
\[\cN_{1,m}\cong \Sym^{m-1}\big(\P^1\big) \cong \P^{m-1}.\]
Moreover, the coef\/f\/icient system $\Phi_{1,m}$ is (up to shift) just the trivial sheaf,
so its (hyper)coho\-mo\-logy
is the cohomology of a smooth projective variety with trivial coef\/f\/icients, carrying in each degree a pure Hodge
structure of the correct weight. It is indeed straightforward to check that in $T$-degree $l=1$, the
series~\eqref{eq_mmns} simply gives the (shifted) weight polynomials~\eqref{eq_Pn}
of projective spaces in each degree in~$q$.
\end{Example}

\begin{Example} \label{ex_N24}
For the case of the curve class having multiplicity two, in the simplest case $\cN_{2,3}\cong\P^1$ is still
nonsingular. The next space $\cN_{2,4}$ is more interesting. The reduced variety underlying~$\cN_{2,4}$
is isomorphic to $\P^3$. As discussed in~\cite[Section~4.1]{PT} however, this cannot be the full answer:
the numerical invariant $P_{2,4}$ equals $4$ and not $-4$, which would be the answer in case
the moduli space were just a smooth $\P^3$. The moduli space $\cN_{2,4}$ is in fact a non-reduced scheme,
the tickening of $\P^3$ along the embedded quadric $Q\subset\P^3$, with Zariski tangent spaces of
dimension $4$ along the quadric. Via torus localization, this indeed gives the correct value
$P_{2,4}=4$. When expanded, the ref\/ined expression, the coef\/f\/icient of $T^2 q^4$ in~\eqref{eq_mmns}, is
$t +2 + t^{-1}$, representing the (appropriately shifted) cohomology of the quadric $Q$ itself. I now
show, under an assumption, how this answer arises in the present framework.

Recall the function $f_{l,m}$,
whose local derivatives cut out $\cN_{l,m}$ inside the smooth ambient spa\-ce~$N_{l,m}$.
Let us make the assumption that locally around every closed point $p\in Q\subset \cN_{2,4}\subset N_{2,4}$ of
the quadric, there are local (analytic) coordinates $x_1, \ldots, x_{2n}$ on the smooth $N_{2,4}$, such that
the function $f_{2,4}$ is locally of the form
$f_{2,4}(x_1, \ldots, x_{2n}) = x_3x_4^2 + \sum\limits_{i=5}^{2n} x_i^2.$
In this case, the degeneracy locus $\cN_{2,4}$ indeed looks around $p$ like a smooth threefold, parametrized
by $x_1$, $x_2$, $x_3$, with
an embedded nilpotent direction with coordinate $x_4$ along the codimension one locus $\{x_3=0\}$.
While it should be possible to prove this assumption
along the lines of~\cite[Section~3.3]{DS}, starting from the quiver model, we will not attempt to do that here.

Let us see, under this assumption, what the vanishing cycle module $\Phi_{2,4}$ looks like. This module
lives on the reduced degeneracy locus $\P^3$. At smooth points, i.e.\ away from the quadric~$Q$,
we must have a rank-one local system. At points of~$Q$, the local structure is described by
Proposition~\ref{prop_xy2} below. Globally, the summand with nontrivial
monodromy can be described as follows. Let $\pi\colon\tilde Q\to \P^3$ be the double cover
of $\P^3$ branched along $Q$. Let
$U=\P^3\setminus Q \stackrel{j}{\hookrightarrow}\P^3$ be the inclusion of the complement of $Q$, and
$\tilde U = \pi^{-1}(U)$. Then over $U$, we have a mixed Hodge module $L$, with underlying local system of rank one
with nontrivial $\Z/2$ monodromy, such that \[(\pi|_{\tilde U})_* \Q^H_{\tilde U} = \Q^H_{\tilde U} \oplus L\] and
thus
\begin{equation}\label{eq_sum_P3}
\pi_* \Q^H_{\tilde Q} = \Q^H_{\P^3} \oplus  j_{!*} L.
\end{equation}
With this notation, we have, as a consequence of Proposition~\ref{prop_xy2}, that
\begin{equation}
\Phi_{2,4} \cong \Q_Q^H(1)[2] \oplus j_{!*}L(2)[2].
\label{eq_phi24}
\end{equation}
In particular, the mixed Hodge module $\Phi_{2,4}$ is a direct sum of pieces with dif\/ferent weights.

On the other hand, $\tilde Q$ is just another quadric, this time in $\P^4$. Moreover
\[\hH^*(\P^3, \pi_* \Q^H_{\tilde Q}) \cong H^*(\tilde Q, \Q)\cong H^*(\P^3, \Q).\]
Thus, from~\eqref{eq_sum_P3}, it follows that
\[\hH^*(\P^3, j_{!*}L)=0.\]
We f\/inally deduce from~\eqref{eq_phi24} that
\[ \hH^*\!\!\left(\P^3, \Phi_{2,4}\right) \cong H^*\!\!\left(Q, \Q^H_{Q}(1)[2]\right).
\]
Thus indeed, the Hodge structure on this cohomology is pure with the correct degrees,
and the weight polynomial is $t+ 2 + t^{-1}$ as demanded by our formula~\eqref{eq_mmns}.
\end{Example}

\begin{Exercise} Study the explicit geometry of the next moduli space $\cN_{2,5}$ and the mixed Hodge module $\Phi_{2,5}$ in similar detail.
\end{Exercise}

\section{Vector spaces underlying partition functions}

The aim in this section is to investigate whether the relationship~\eqref{fullZinverse}
can be used to study the vector spaces underlying these partition functions. As discussed before,
on the left hand side of~\eqref{fullZinverse}, the four-dimensional partition
function $Z_{\C^2, r=1}(q_1, q_2,\Lambda)$ is the graded character of an actual triply-graded vector space
\begin{equation} \displaystyle\bigoplus_{k\geq 0} H^0(\O_{S^k(\C^2)}) \cong S^*\C[x,y].
\label{iso:2d}
\end{equation}
On the right hand side of~\eqref{fullZinverse}, we have weight polynomials.
In general, forming the weight polynomial involves taking an Euler characteristic: in the def\/ining
formula~\eqref{Wdef}, there can be cancellation between dif\/ferent pieces of cohomology with the same weight
(see Examples~\ref{ex:ellcurve} and~\ref{ex:A2ish}).
However, it is clear from the formula that this cannot occur for pure Hodge structures, where the weight
f\/iltration is trivial, agreeing with the obvious grading of cohomology by degree.
Indeed, by Corollary~\ref{cor:pure}(1), all these cohomologies carry pure Hodge structures, and thus
\[ Z_{X}(q,T,t) = \cha_{(\C^*)^3} \bigoplus_{l,m} \hH^*\!\left(\cN_{l,m}, \Phi_{l,m}\right)\]
is also the graded character of a triply-graded (super) vector space, graded by the curve degree~$l$,
the point degree~$m$ and the cohomology degree.

Given that the four-dimensional partition function is the graded character of a symmetric space of a vector space,
the inverse operation in~\eqref{fullZinverse} has a natural interpretation: we obtain
\[ Z_{X}(q,T,t) = \cha_{(\C^*)^3}  \Lambda^*\C[x,y]
\]
under a grading where $x$, $y$ are viewed as {\rm odd} variables (introducing signs into the Hilbert series)
of weights $(1,0,\pm\half)$, and the outer exterior operation
has weight $(1,1,0)$. The following is the main result of our paper.
\begin{Theorem} \label{thm:iso} There exists a natural $\GL(2)\times\C^*$-equivariant isomorphism
\begin{equation} \bigoplus_{l,m} \hH^*\!\left(\cN_{l,m}, \Phi_{l,m}\right) \cong \Lambda^*\C[x,y]
\label{iso:3d}
\end{equation}
of $($super$)$ vector spaces.
\end{Theorem}
\begin{proof} Our computation of the torus characters means that we certainly have a
triply graded, in other words $(\C^*)^3$-equivarient isomorphism.

Return to the four-dimensional partition function for a moment. The $T=(\C^*)^2$-action on~$\C^2$ is part of a $\GL(2)$-action;
the element $\left(\begin{matrix} 0 & 1 \\ 1 & 0\end{matrix}\right)$ acts on the partition function simply by
interchanging $x$ and $y$. Looking at the change of variables in Corollary~\ref{cor_identify}, we see that
this corresponds to mapping $t^{\half}\mapsto t^{-\half}$, while keeping the other variables f\/ixed,
in other words to Poincar\'e duality on
the individual cohomologies $\hH^* (\cN_{l,m}, \Phi_{l,m})$, for every f\/ixed~$l$,~$m$. This
then means that we can enhance the corresponding $(\C^*)^2$-action on the left hand side
of~\eqref{iso:3d} to a~$\GL(2)$-action, integrating the $\sl(2)$-action coming from
Corollary~\ref{cor:pure}(2), with Weyl element given by the action of the Hard Lefschetz isomorphism. We
then have compatible $\GL(2)\times\C^*$-actions on all the spaces in~\eqref{iso:2d},~\eqref{iso:3d}, compatible
with the isomorphisms.
\end{proof}

Before I discuss~\eqref{iso:3d} any further, let me make a digression. Recall Theorem~\ref{thm_mmns:theta},
computing a six-dimensional partition function depending on a parameter $\theta$. Checking the ef\/fect of
the truncation, the following extension of Theorem~\ref{thm:iso} turns out to be
compatible with the partition function.

\begin{Theorem} \label{thm:iso:theta} For $n\geq 0$ and $\theta\in\cC_n$, there exists a $\GL(2)\times\C^*$-equivarient isomorphism
\begin{equation} \bigoplus_{l,m} \hH^*\big(\cN^\theta_{l,m}, \Phi^{\theta}_{l,m}\big) \cong \Lambda^*\C[x,y]_{n-1}
\label{iso:3d:theta} \end{equation}
of $($super$)$ vector spaces, where $\C[x,y]_n$ denotes the space of polynomials of total degree at most~$n$.
\end{Theorem}

I want to address the issue of naturality of the isomorphisms~\eqref{iso:3d},~\eqref{iso:3d:theta}.
Recall one f\/inal time the conifold quiver with potential $(Q,W)$ and its close relative, the framed
quiver $(\tilde Q, \tilde W)$. In a recent paper~\cite{KS2}, Kontsevich and Soibelman introduce
an associative algebra $\cH(Q,W)$, the (critical) {\em cohomological Hall algebra} of the quiver $(Q,W)$, built
from {\em all} representations of the quiver~$Q$, the vanishing cycle complex of
the function~$\Tr(W)$ def\/ined on these spaces of representations, and equivariant cohomology with respect to
the natural action of products of general linear groups.

The associative product on $\cH(Q,W)$ is def\/ined
in~\cite{KS2} by using a version of the standard diagram which essentially fuses two representations $\cR_i$ of $(Q,W)$
into a third one: an exact sequence
\[ 0\to \cR_2 \to \cR_3 \to \cR_1 \to 0\]
of representations of $(Q,W)$ gives the multiplication
\begin{figure}[h]
\centering
\includegraphics{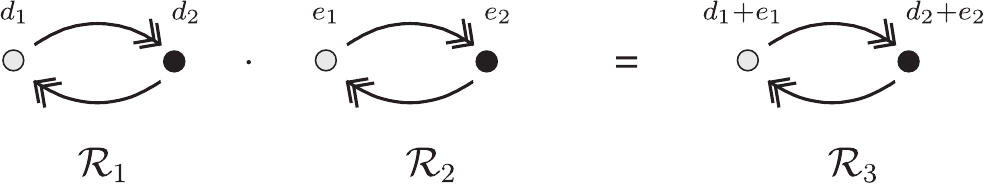}
\end{figure}

The same construction should turn the space on the right hand side of~\eqref{iso:3d},
the sum of cohomologies of representation spaces of the quiver $(\tilde Q, \tilde W)$, into {\em modules} over the
algebra $\cH(Q,W)$. The point is that the dimension on the extra vertex in $\tilde Q$ is always one, and
this is unchanged by the operations which attach representations $\cR$ of $(Q,W)$: an exact sequence
\[ 0\to \cR \to \tilde\cR_2 \to \tilde\cR_1 \to 0
\]
of representations of $(\tilde Q,\tilde W)$ should give the multiplication
\begin{figure}[h]
\centering
\includegraphics{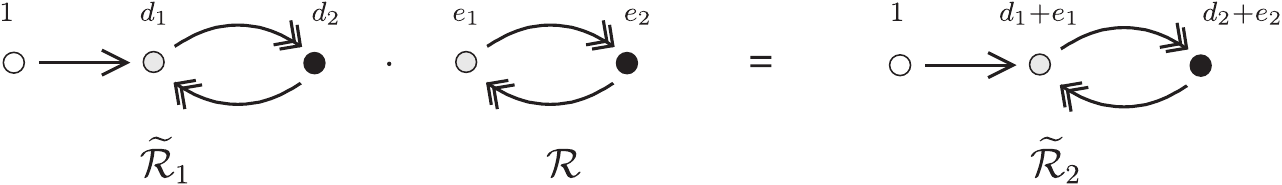}
\end{figure}

\noindent
See the recent preprint~\cite[Section~2]{GS} for a similar, possibly not unrelated, perhaps mirror, idea of how to construct
representations of the conifold CoHA (algebra of BPS states).

The structure of the algebra $\cH(Q,W)$ is not known explicitly at present. But the most natural
extension of the ideas above is that there should exist unique isomorphisms~\eqref{iso:3d},~\eqref{iso:3d:theta}
of $\cH(Q,W)$-modules, exhibiting $\Lambda^*\C[x,y]$ as a Verma-type module of the algebra $\cH(Q,W)$ and
$\Lambda^*\C[x,y]_n$ as f\/inite-dimensional highest-weight quotients.

\begin{Remark} As a f\/inal point, notice that the fact that we could interpret both the four- and the
six-dimensional partition function as graded characters of vector spaces (as opposed to virtual torus representations)
depended on {\em both} sides on a fortuitous lack of cancellation. On the six-dimensional side, this is
provided by the purity of cohomologies of
mixed Hodge modules of vanishing cycles. On the four-dimensional side,
\eqref{eq_nek_def} involves an index, in other words an alternating sum of cohomologies; however,
as discussed around~\eqref{eq_nek_H0}, there is no higher cohomology, leading to an actual, rather than virtual,
torus representation.
I can see no direct connection between these facts, but it would be really interesting if one existed.
\end{Remark}

\appendix

\section{Weights on cohomology}
\label{app_weights}

Since the weight f\/iltration plays an important role in the paper, here we collect some conventions for weights.
Recall that by Deligne's theorem, the rational cohomology $H^*(X, \Q)$ of a~variety~$X$, not necessarily smooth
or projective, admits a mixed Hodge structure: it carries a weight f\/iltration, so that the associated graded pieces carry
pure Hodge structures. If $X$ is smooth and projective, then the weight f\/iltration is trivial: it coincides with the
degree f\/iltration. More generally, the smoothness or projectivity on $X$ imply estimates on the weights appearing
in the cohomology.

The modern way to view Deligne's mixed Hodge structure is to f\/irst consider, for a ge\-ne\-ral scheme~$X$,
the category of pure Hodge modules~$\HM(X)$, sitting
inside a category of mixed Hodge modules $\MHM(X)$; objects in these categories are complicated, represented by a~perverse
topological $\Q$-sheaf,
an associated complex algebraic gadget (a $D$-module) and some compatibility data. Pure Hodge modules are direct
sums of Hodge modules, each of which
has a f\/ixed weight; mixed Hodge modules are extensions of pure Hodge modules. The next step is to form complexes
of mixed Hodge modules, to arrive at the bounded derived category $\D^b\MHM(X)$. For an object $\F\in\D^b\MHM(X)$,
$\F[k]$ denotes the complex~$\F$ shifted by $k$~places.

Maps $f\colon X\to Y$
between algebraic varieties def\/ine pushforward maps between the corresponding derived categories; there are two types of
pushforward $f_*$ and $f_!$, corresponding in the case the target $Y=\{P\}$ is a point to taking cohomology on~$X$ and
cohomology with compact support. The categories of Hodge modules and mixed Hodge modules over a point are equivalent to the
categories of pure and mixed Hodge structures; thus cohomologies of mixed Hodge modules on arbitrary~$X$
carry mixed Hodge structures. There are also two types of pullback maps~$f^*$ and~$f^!$.

Every smooth variety $X$ carries a canonical pure Hodge module $\Q_X^H\in\HM(X)$ of weight $2\dim X$. Proper pushforward
preserves purity, so we obtain that if $X$ is moreover proper, then $H^*(X,\Q_X^H)$ is pure; this is simply the classical
cohomology of~$X$. For the projective line, we have
\[ H^*(\P^1, \Q_{\P^1}^H)\cong \Q \oplus \Q(-1)[-2];
\]
here $\Q(-1)$ is the Tate Hodge structure, a one-dimensional pure Hodge structure of weight $2$ (the analogue of the motive
$\L$ in the motivic ring of varieties). The shift $[-2]$ corresponds to the fact that it arises as the second cohomology.
For an arbitrary complex of mixed Hodge modules $\F\in\D^b\MHM(X)$, we let $\F(i)$ denote the tensor product
of $\F$ with $f^*\Q(i)$, where $f\colon X\to\{P\}$ is the structure morphism.
If $\F$ is pure of weight $k$, then $\F(i)$ is pure of weight $k-2i$.

To treat odd-dimensional varieties, it is convenient to adjoin to the category of mixed Hodge structures a half-Tate
object $\Q(\half)$, of weight $-1$, with the property that $\Q(\half)\otimes\Q(\half)\cong\Q(1)$;
see~\cite[Section~3.4]{KS2}. Then for any smooth projective $X$, the cohomology
\[H^*(X, \Q^H_X(\dim X/2)[\dim X])\] lives in palindromic degrees $-\dim X, \ldots, \dim X$ with the same weights.
The pullback to any variety $X$ of $\Q(\half)$ under the structure morphism will still be denoted by $\Q(\half)$;
we denote by $\tMHM(X)$ the extended category of mixed Hodge modules on a variety~$X$.

Given a complex of mixed Hodge modules $\F\in\D^b\tMHM(X)$, consider the compactly supported
(hyper)cohomology $\hH_c^*(X, \F)$ with its weight f\/iltration $W_n$.
We can consider the weight polynomial (sometimes called Serre polynomial)
\begin{equation*}
W\big(X, \F; t^\half\big) = \sum_{n\in \Z} t^{\frac{n}{2}} \sum_{i\in\Z} (-1)^i \dim \Gr^W_{n} \hH_c^i(X, \F)\in\Z\big[t^{\pm\half}\big].
\end{equation*}
Then
\[W\big(X, \F(i)[j]; t^\half\big)= (-1)^j t^{-\frac{i}{2}} W\big(X, \F; t^\half\big).
\]
A very important property of the weight polynomial is its additivity under stratif\/ications:
because of the long exact sequence in cohomology
and its compatibility with Hodge and weight f\/iltrations, the weight polynomial behaves well
under decompositions $X= (X\setminus Z) \cup Z$, where $Z\subset X$ is a closed subvariety.

\begin{Example} To start with,
\[
W\big(\A^n, \Q_{\A^1}^H; t^\half\big) = t^n
\]
as the only nonzero compactly supported cohomology of ${\mathbb A}^n$ is one-dimensional of weight $2n$ in degree $2n$.
Also
\[
W\big(\P^1, \Q_{\P^1}^H; t^\half\big) = 1 + t
\]
and so
\[
W\big(\P^1, \Q_{\P^1}^H(1/2)[1]; t^\half\big) = -\big(t^{\half} + t^{\half}\big);
\]
more generally,
\begin{equation}\label{eq_Pn}
W\big(\P^n, \Q_{\P^n}^H(n/2)[n]; t^\half\big)
= (-1)^n \big(t^{-\frac{n}{2}} + \cdots + t^{\frac{n}{2}} \big) = (-1)^n  \frac{t^{\frac{n+1}{2}}  - t^{-\frac{n+1}{2}}}{t^{\frac{1}{2}}  - t^{-\frac{1}{2}}}.
\end{equation}
Even more generally, for a smooth projective variety $X$, $H^i(X, \Q)$ is of weight $i$, so we have
\[W\big(X, \Q_{X}^H; t^\half\big) = \sum_{n=0}^{2\dim X} (-1)^i b_i(X) t^{\frac{i}{2}},
\]
where the $b_i$ are the Betti numbers of $X$.
\end{Example}
\begin{Example}\label{ex:ellcurve} For $E$ an elliptic curve,
\[W\big(E, \Q_{E}^H; t^\half\big) = 1- 2 t^\half + t.
\]
Consider now the complex of Hodge structures (complex of Hodge modules on a point $P$)
$\F=H^*(E, \Q)\oplus \Q(1/2)^{\oplus 2}$. This has the
previous pieces $H^{2i}(E, \Q)[-2i]$ in weight $2i$, and $H^1(E)[-1]\oplus \Q^{\oplus 2}$ in weight $1$.
Thus, its weight polynomial is
\[W\big(P, \F; t^\half\big) = \big(1- 2 t^\half + t\big) + 2 t^\half = 1 + t,
\]
indistinguishable from the weight polynomial coming from the cohomology of $\P^1$ with constant coef\/f\/icients.
This is the kind of cancellation that cannot occur under the purity statement of Theorem~\ref{thm:pure}.
\end{Example}

\begin{Exercise}\label{ex:A2ish} For a more geometric example of cancellation,
let $X$ be the blowup of a point in $\A^1\times\C^*$. Show that its
weight polynomial with constant coef\/f\/icients is
\[ W\big(X, \Q_{X}^H; t^\half\big) = t^2,
\]
indistinguishable from that of $\A^2$, even though $X$ obviously has odd cohomology.
\end{Exercise}

\section{A vanishing cycle computation}
We compute the vanishing cycle mixed Hodge module of the function $f(x,y)=xy^2$ on $X=\C^2$.
This sort of computation is presumably trivial for the experts, but I couldn't f\/ind the answer in the literature
in a form needed above.

Start with the inclusions
\[P=\{0\}\in Z=\{df=0\}_{\rm red}\subset Y =\{f=0\}_{\rm red}\subset X,\]
with $P$ being the origin, $Z$ the $x$-axis and $Y$ the union of the axes.
On the smooth $X=\C^2$, we have $\Q_X^H[2]\in\MHM(X)$, a pure Hodge module of weight $2$.
On the hypersurface $Y$, we have the mixed Hodge module of vanishing cycles
$\phi_f\Q_X^H\in\MHM(Y)$. There is an exact sequence of MHM's on~$Y$
\[0 \to \Q^H_Y[1] \to \psi_f\Q_X^H[2] \to \phi_f\Q_X^H[2] \to 0.\]
The weight f\/iltration of $\Q^H_Y[2]$ has $IC^H_Y$ in weight~1, then $i_*\Q_P^H$ in weight~0, with
$i\colon P \to Y$ denoting the inclusion of the origin~$P$.
Now this injects into $\psi_f\Q_X^H[2]$, and the action of the monodromy operator~$N$ on
$\psi_f\Q_X^H[2]$ has to satisfy Hard Lefschetz. The nearby cycle over the
punctured $x$-axis consists of two points.
Correspondingly, there is a local system $M$ on the smooth locus $j\colon U\to Y$,
which is trivial over the punctured $y$-axis, and is of rank two over the punctured $x$-axis. Then
$\psi_f\Q_X^H[3]$ has to look as follows:
in weight~$2$: $i_*\Q_P^H(-1)$; in weight~$1$: $IC^H_Y(M)$; in weight $0$: $i_*\Q_P^H$. We have $N^2=0$, but $N$
is nonzero, taking the weight-$2$ piece isomorphically to the weight-$0$ piece.
The cokernel of the inclusion $\Q^H_Y[2] \to \psi_f\Q_X^H[3]$ is $\phi_f\Q_X^H[3]$. We deduce
that $\phi_f\Q_X^H[3]$ has weight f\/iltration
as follows. In weight~$2$, the associated graded is $i_*\Q_P^H(-1)$; in weight~$1$, it is $IC^H_Z(L)$, where~$L$ is the
nontrivial rank-one local system on the punctured $x$-axis with $\Z/2$ monodromy.
Note f\/inally that the weight f\/iltration splits, either because of duality, or due to the fact that
the monodromy-invariant part $\phi_{f,1}\Q_X^H[2] \cong i_*\Q^H(-1)$ is always a direct summand. In summary,
\begin{Proposition} \label{prop_xy2}
In the category $\MHM(Z)$, the module $\phi_f\Q_X^H[2]$ splits as a direct sum
\[
\phi_f\Q_X^H[2] \cong i_*\Q^H_P(-1) \oplus IC^H_Z(L),
\]
with weights $2$ and $1$ respectively.
\end{Proposition}

\subsection*{Acknowledgements}

I~wish to thank Jim Bryan, Lotte Hollands, Dominic Joyce, Davesh Maulik,
Geordie Williamson and especially Ian Grojnowski for comments and discussions. This research was supported by
EPSRC Programme Grant EP/I033343/1, and by a Fellowship from the Alexander von Humboldt Foundation.
Part of this paper was prepared while I was visiting the Department of Mathematics, Freie Universit\"at
Berlin; I wish to thank them and especially Klaus Altmann for hospitality.

\pdfbookmark[1]{References}{ref}
\LastPageEnding

\end{document}